\documentclass{amsart}
\usepackage{latexsym}
\usepackage{amsmath, amssymb, amsfonts, amsxtra}
\usepackage{amsthm}

\author{Richard Oberlin}
\email{oberlin@math.wisc.edu}
\address{University of Wisconsin-Madison, Mathematics Department, 480 Lincoln Dr, Madison WI 53706}

\title{Two bounds for the x-ray transform}
\date{}

\newcommand {\leb}{\mathcal{L}}

\newcommand {\rea}{\mathbb{R}}
\newcommand {\sph}{\mathbb{S}}
\newcommand {\ma}{\mathcal{M}}
\newcommand {\proj}{\mathrm{proj}}
\newcommand {\spa}{\mathop{\mathrm{span}}}
\newcommand {\composed} {\circ}

\newcommand {\mass}{\mathfrak{M}}

\newtheorem{theorem}{Theorem}[section]
\newtheorem{corollary}{Corollary}[section]
\newtheorem{proposition}{Proposition}[section]
\newtheorem{claim}{Claim}[section]

\begin{document}

\begin{abstract}
We use the arithmetic-combinatorial method of Katz and Tao to give
 mixed-norm estimates for the $x$-ray transform on $\rea^d$ 
when $d \geq 4.$ As an application, we obtain an improved estimate 
for the Hausdorff dimension of $(d,k)$ sets, which are subsets of $\rea^d$
containing a translate of every $k$-plane. 
\end{abstract}

\maketitle

\section{Introduction}

Let $e_1,\ldots,e_d$ be an orthonormal basis for $\rea^d$ and let $H = \spa(e_1,\ldots,e_{d-1})$.
For $\xi, x \in H$ and a function $f$ defined on $\rea^d$, 
the local $x$-ray transform of $f$ at the line $x + \rea(\xi + e_d)$ may be defined
\[
T[f](\xi,x) = \int_0^1 f(x + t(\xi + e_d))\ dt.
\]
We consider the ``Kakeya-order'' 
mixed norm of $T[f]$
\begin{equation} \label{kmndef}
\|T[f]\|_{L^q(L^r),K} = \left(\int_{B(0,C)}\left(\int_{\rea^{d-1}} |T[f](\xi,x)|^r\ dx\right)^{\frac{q}{r}}\ d\xi\right)^{\frac{1}{q}}
\end{equation}
where $C > 0$ and $B(0,C)$ denotes the ball centered at $0$ with radius $C$ in $\rea^{d-1}$, and  
we aim to prove bounds of the form
\begin{equation} \label{kakeyaorder}
\|T[f]\|_{L^q(L^r),K} 
\lesssim \|f\|_{L^p(\rea^d)}.
\end{equation}
Through a covering argument, one may see that the bound (\ref{kakeyaorder}) is equivalent to the more conventional version
\begin{equation} \label{kakeyaorderconv}
\left(\int_{\sph^{d-1}} \left(\int_{\xi^\perp} \left|\int_{\rea} f(x + t\xi)\ dt\right|^r dx \right)^{\frac{q}{r}} d\Omega(\xi)\right)^{\frac{1}{q}} \lesssim \|f\|_{L^{p}(\rea^d)} 
\end{equation} 
where $\sph^{d-1}$ is the $d-1$ sphere, $\Omega$ is surface measure on $\sph^{d-1}$, and we only consider $f$ supported on a fixed ball.

By testing $T$ on the characteristic functions of $\delta$-neighborhoods of points and line segments, and letting $\delta$ approach $0$, one sees that 
\begin{equation} \label{prcondition}
1 + \frac{d-1}{r} \geq \frac{d}{p}
\end{equation} 
and
\begin{equation} \label{pqrcondition}
\frac{1}{q} + \frac{1}{r} \geq \frac{1}{p}
\end{equation}
are necessary conditions for the bound (\ref{kakeyaorder}) to hold. It is 
conjectured that, together with the condition $r < \infty$, these are 
also sufficient. This was shown to be the case for $p < \frac{d+1}{2}$ by Drury in \cite{dr} and for $p = \frac{d+1}{2}$ by Christ in \cite{ch}. 

For $\delta > 0, \xi \in \sph^{d-1}$, and $a \in \rea^d$, let $\tau_{\delta}(\xi,a)$ denote the $\delta$-neighborhood of the line segment with endpoints $a,a+\xi$. Then the Kakeya maximal operator is 
defined 
\[
K[f](\xi) = \sup_{a \in \rea^d} \frac{1}{\leb^d(\tau_{\delta}(\xi,a))} \int_{\tau_{\delta}(\xi,a)} f(y)\ dy
\]
where we use $\leb^{d}$ to denote Lebesgue measure in $\rea^d$. 
Through an application of H\"{o}lder's inequality to the inner norm, one may see that the bound (\ref{kakeyaorderconv}) implies the bound for the Kakeya maximal operator
\begin{equation} \label{kakeyamobound}
\|K[f]\|_{L^q(\sph^{d-1})} \lesssim \delta^{-\frac{d-1}{r}} \|f\|_{L^p(\rea^d)}.
\end{equation}
This implication does not seem to be reversible. However, it is reasonable 
to expect that strategies employed for proving (\ref{kakeyamobound}) may 
extend to yield (\ref{kakeyaorder}). For example, in \cite{wo2} Wolff showed that (\ref{kakeyamobound}) holds with $(p,q,r) = (p_w,q_w,r_w - \epsilon),$ where
\[
p_w = \frac{d+2}{2},\ q_w = \frac{(d-1)(d+2)}{d},\ r_w = \frac{(d-1)(d+2)}{d-2} ,
\] 
and later, in \cite{wo}, he refined the technique to show that when $d=3$
\[
\|T[f]\|_{L^{q_w}(L^{r_w}),K} \lesssim \|f\|_{p_w,\epsilon}
\]
where $\|\cdot\|_{p,\epsilon}$ denotes the inhomogeneous Sobolev norm with $\epsilon$ derivatives in $L^p$. In \cite{lt}, {\L}aba and Tao extended this
result to higher dimensions, obtaining for $d \geq 4$ 
\[
\|T[f]\|_{L^{q_w}(L^r),K} \lesssim \|f\|_{p_w,\alpha}
\]
when $r = 2(d+2),\ \alpha = \frac{d-3}{2(d+2)} + \epsilon.$

In \cite{kt}, Katz and Tao showed that (\ref{kakeyamobound}) holds with $p = \frac{4d+3}{7}, q = \frac{4d+3}{4},$ and $r = \frac{4d+3}{3} - \epsilon$. We answer
the question of whether this may be extended to a result of type (\ref{kakeyaorder}) affirmatively.
\begin{theorem} \label{maintheorem}
When $d \geq 6$ and $\epsilon > 0$, the bound (\ref{kakeyaorder}) holds with $p = \frac{4d+3}{7}, q = \frac{4d+3}{4}- \epsilon,$ and $r = \frac{4d+3}{3} - \epsilon$.
\end{theorem}
Except for the losses of $\epsilon$, Theorem \ref{maintheorem} is optimal for the ratio $\frac{r}{p} = \frac{7}{3}$ in the sense that it obtains the smallest $p$ and the largest $q$ permitted by (\ref{prcondition}) and (\ref{pqrcondition}).
In contrast to \cite{wo} and \cite{lt}, the 
refinement of the argument in \cite{kt} needed to obtain Theorem \ref{maintheorem} is quite minimal.  
 
A measurable subset of $\rea^{d}$ is said to be a Kakeya set if it contains
a unit line segment in every direction. Bourgain showed in \cite{bg} that
the bound (\ref{kakeyamobound}) implies that Kakeya sets must have Hausdorff dimension at least $d - (d-1)\frac{p}{r}.$ In \cite{kt}, Katz and Tao
showed that the Hausdorff dimension of a Kakeya set must be at least $(2 - \sqrt{2})(d-4) + 3,$ which is stronger than the result implied by their maximal operator bound. To obtain this dimension estimate, they combined the main estimate from their maximal operator bound with $2$-dimensional Kakeya methods, as in the ``hairbrush'' construction of \cite{wo2}, and used an iterative technique to improve the resulting estimate. We iterate the main estimate from Theorem \ref{maintheorem} to obtain   
\begin{theorem} \label{iteratedtheorem}
For $d \geq 4$ and $\epsilon > 0$, there exist $p_\epsilon, q_\epsilon, r_\epsilon$ so that (\ref{kakeyaorder}) holds with $(p,q,r) = (p_\epsilon,q_\epsilon,r_\epsilon)$ where
\[
\frac{r_\epsilon}{p_\epsilon} > 1 + \sqrt{2} - \epsilon, \text{\ and\ } \frac{q_\epsilon}{p_\epsilon} > 1 + \frac{\sqrt{2}}{2} - \epsilon.
\]
\end{theorem}
The value of $p_\epsilon$ required in Theorem \ref{iteratedtheorem} is nonoptimal for
the ratio $\frac{r}{p} = 1 + \sqrt{2}$, and in fact approaches $\infty$ as $\epsilon$ approaches $0$.

Let $G(d,k)$ denote the set of $k$-dimensional subspaces of $\rea^d$. Then, for $M \in G(d,k)$ and $x \in M^{\perp}$, 
the $k$-plane transform of $f$ at $M+x$ is defined 
\[
T^k[f](M+x) = \int_{M + x} f(y)\ dy.
\]
The local $k$-plane transform conjecture for $k > 1$ is that when
\[
k + \frac{d-k}{r} \geq \frac{d}{p}, \text{\ \ \ and\ \ \ } \frac{k}{q} + \frac{1}{r} \geq \frac{1}{p},
\]
the bound 
\begin{equation} \label{kplanekakeyabound}
\left(\int_{G(d,k)} \left( \int_{M^{\perp}} |T^k[f](M + x)|^r \ dx \right)^{\frac{q}{r}}\ dM \right)^{\frac{1}{q}} \lesssim \|f\|_{L^p(\rea^d)}
\end{equation}
holds for $f$ supported on a ball, where we integrate with respect to rotation invariant measure on $G(d,k)$.
Christ showed in \cite{ch} that this holds for $p \leq \frac{d+1}{k+1}$.

In certain cases, we may improve this range of $p$ by applying the method of \cite{ro} to Theorem \ref{maintheorem}. 
\begin{corollary} \label{kplanecorollary}
Suppose $k \geq 2$, $d \geq 5 + k$, and 
\[
p = \frac{d}{k+\frac{3}{4}},\ 1 + \frac{d-k}{r} > \frac{d}{p},\ q < \frac{4(d-(k-1))+3}{4}.
\]
Then (\ref{kplanekakeyabound}) holds for $f$ supported on a ball.
\end{corollary}
A limitation of the method in \cite{ro}, is that we are not able to obtain 
the full range of $q$ above. Also, the value of $p$ above is approximate and may be slightly improved, although we omit the numerological details.

By a similar method, Theorem \ref{iteratedtheorem} yields
\begin{corollary} \label{unsharpkplanecorollary}
Let $\epsilon > 0$, $k \geq 2$, $d \geq 3+k$. Then, there exist $p_\epsilon,r_\epsilon$ satisfying
\[
r_\epsilon > p_\epsilon (1 + \sqrt{2} - \epsilon)^k
\]
so that (\ref{kplanekakeyabound}) holds for, say, $(p,q,r) = (p_\epsilon,1,r_\epsilon)$ and $f$ supported on a ball.   
\end{corollary}   
A $(d,k)$ set is a measurable subset of $\rea^d$ which contains a translate of every $k$-plane. Our main interest in Corollary \ref{unsharpkplanecorollary} is that it yields
\begin{corollary} \label{dksetcorollary}
The Hausdorff dimension of every $(d,k)$ set is at least 
\[
d - \frac{d-k}{(1 + \sqrt{2})^k}.
\]
Furthermore, $(d,k)$ sets have positive measure if  
\[
(1 + \sqrt{2})^{k-1} + k > d.
\]
\end{corollary}
This may be seen as a natural update of Bourgain's result from \cite{bg} that 
$(d,k)$ sets have positive measure for $2^{k-1} + k \geq d.$

Instead of the Kayeya-order mixed norms (\ref{kmndef}), one may consider the Nikodym-order mixed norms of $T[f]$
\[
\|T[f]\|_{L^q(L^r),N} = \left(\int_{\rea^{d-1}}\left(\int_{B(0,C)} |T[f](\xi,x)|^r\ d\xi\right)^{\frac{q}{r}}\ dx\right)^{\frac{1}{q}}.
\]
In order for the bound
\begin{equation} \label{nikodymorder}
\|T[f]\|_{L^q(L^r),N} 
\lesssim \|f\|_{L^p(\rea^d)}
\end{equation}
to hold, we again have the necessary conditions (\ref{prcondition}), (\ref{pqrcondition}), and $r < \infty$. Unless we impose the additional assumption that $f$ is supported away from $H$, we have another necessary condition
\begin{equation} \label{qpcondition}
1 + \frac{d-1}{q} \geq \frac{d}{p}
\end{equation}
which follows from the application of $T$ to characteristic functions of $\delta$-neighborhoods of a point in $H$. Tao showed in \cite{tao} that bounds for the Kakeya and Nikodym maximal operators are roughly equivalent. We 
observe that his proof carries over to the general mixed norm case, and hence combined with Theorem \ref{maintheorem} yields
\begin{corollary} \label{nikodymcorollary}
When $d \geq 6$ and $\epsilon > 0$, the bound (\ref{nikodymorder}) holds with $p = \frac{4d+3}{7}, q = \frac{4d+3}{4}- \epsilon,$ and $r = \frac{4d+3}{3} - \epsilon$.
\end{corollary}
One may also formulate a corresponding version of Theorem \ref{iteratedtheorem}.

We prove Theorem \ref{maintheorem} in sections \ref{rwesection}, \ref{twoendssection}, \ref{continslice}, and \ref{pigeonholesection}.
We give the additional arguments needed for Theorem \ref{iteratedtheorem} in section \ref{pigeonholesection2}. We show that Corollary \ref{nikodymcorollary} follows from Theorem \ref{maintheorem} in section \ref{projectivesection}. Corollaries \ref{kplanecorollary}, \ref{unsharpkplanecorollary}, and \ref{dksetcorollary} follow from the arguments in \cite{ro}.

\section{Reduction to weak estimates} \label{rwesection}
We first note that, since $T$ is local, when $p \leq q \leq r$ it suffices to prove (\ref{kakeyaorder}) for $f$ supported on the 
cube $Q$ centered at $\frac{1}{2}e_d$ with side length $1$. A natural simplification of the bound (\ref{kakeyaorder})
is the estimate
\begin{equation} \label{rwtestimate1}
\|\lambda \chi_F\|_{L^q(L^r),K}
\lesssim \|\chi_{E}\|_{L^p(\rea^d)}
\end{equation}
where $E \subset Q$, $0 < \lambda \leq 1$, and  
\begin{equation} \label{EFlambdahyp} 
T[\chi_E](\xi,x) \geq \lambda \text{\ for\ } (\xi,x) \in F \subset B(0,C) \times \rea^{d-1}.
\end{equation}
We further simplify the estimate (\ref{rwtestimate1}) to 
\begin{eqnarray} \label{rwtestimate}
\leb^{d}(E)^\frac{1}{p} &\gtrsim& \lambda \left(\frac{\leb^{2(d-1)}(F)}{\ma(F)}\right)^{\frac{1}{q}}\ma(F)^{\frac{1}{r}} \\ \nonumber
&=&\lambda 
\leb^{2(d-1)}(F)^{\frac{1}{q}}\ma(F)^{\frac{1}{r} - \frac{1}{q}}
\end{eqnarray}
where 
\[
\ma(F) = \sup_{\xi \in B(0,C)} \leb^{d-1}(\{x : (\xi,x) \in F\}).
\]
The crude interpolation argument below shows that (\ref{kakeyaorder}) follows from (\ref{rwtestimate}).

\begin{claim} \label{rweclaim}
Suppose that $r > q$ and that the estimate (\ref{rwtestimate}) holds for all $E$ contained in $Q$. Then for any $\epsilon > 0$,
the bound
\[
\|T[f]\|_{L^{q}(L^{r}),K} \lesssim \|f\|_{L^{p+\epsilon}(\rea^d)} 
\]
holds for functions $f$ supported on $Q$.
\end{claim}

\begin{proof}
Without loss of generality, assume that $f$ is nonnegative and 
\begin{equation} \label{normfone}
\|f\|_{L^{p+\epsilon}(\rea^d)} = 1.
\end{equation}
For integers $j,k,l$ let
\begin{gather*}
E_j = \{x \in \rea^d : 2^{j-1} < f(x) \leq 2^j \}
\\ F_{j,k} = \{(\xi,x) \in B(0,C) \times \rea^{d-1} : 2^{k-1} < T[\chi_{E_j}](\xi,x) \leq 2^k \}
\\ F_{j,k,l} = \{(\xi,x) \in F_{j,k} : 2^{l-1} < \leb^{d-1}(\{x' : (\xi,x') \in F_{j,k} \}) \leq 2^l \}.
\end{gather*}
Since each $T[\chi_{E_j}] \leq 1$, we have $F_{j,k} = \emptyset$ for $k > 1$.
Since we only consider $\xi \in B(0,C)$ and $f$ supported on $Q$, $F_{j,k,l} = \emptyset$ for $l > \tilde{C}$. It follows that
\[
\|T[f]\|_{L^q(L^r),K} \lesssim \sum_{j=-\infty}^\infty \sum_{k=-\infty}^0 \sum_{l = -\infty}^{\tilde{C}}
2^{j + k + l(\frac{1}{r} - \frac{1}{q})} \leb^{2(d-1)}(F_{j,k,l})^{\frac{1}{q}}.
\]

Let $\epsilon_1,\epsilon_2 > 0$ and $S_{k,l} = 2^{k \epsilon_1 - l(\frac{1}{r} - \frac{1}{q}) \epsilon_2}$. Then
\begin{align*}
2^{k + l(\frac{1}{r} - \frac{1}{q})} \leb^{2(d-1)}(F_{j,k,l})^{\frac{1}{q}}
&= S_{k,l} 
\left(2^k 2^{l(\frac{1}{r} - \frac{1}{q}) \frac{1 + \epsilon_2}{1 -\epsilon_1}}
\leb^{2(d-1)}(F_{j,k,l})^{\frac{1}{q(1-\epsilon_1)}}\right)^{1-\epsilon_1}.
\end{align*}
Provided that $\epsilon_2$ is sufficiently small relative to $\epsilon_1$, 
we have
\begin{align} \label{rqepsilon2ineq}
\left(\frac{1}{r} - \frac{1}{q}\right)\frac{1+\epsilon_2}{1-\epsilon_1} 
&= \frac{1}{r} - \frac{1}{q(1-\epsilon_1)} + \frac{\epsilon_1}{r(1-\epsilon_1)}
+ \epsilon_2\left(\frac{1}{r(1-\epsilon_1)} - \frac{1}{q(1-\epsilon_1)}\right)
\\ \notag &\geq \frac{1}{r} - \frac{1}{q(1-\epsilon_1)}.    
\end{align}
We then have
\begin{align*}
2^k 2^{l(\frac{1}{r} - \frac{1}{q}) \frac{1 + \epsilon_2}{1 -\epsilon_1}}
\leb^{2(d-1)}(F_{j,k,l})^{\frac{1}{q(1-\epsilon_1)}}
&\lesssim 2^k 2^{l(\frac{1}{r} - \frac{1}{q(1-\epsilon_1)})}
\leb^{2(d-1)}(F_{j,k,l})^{\frac{1}{q(1-\epsilon_1)}}
\\ &\lesssim \leb^{d}(E_j)^{\frac{1}{p}},
\end{align*}
where the first inequality follows from (\ref{rqepsilon2ineq}) and the second inequality follows from (\ref{rwtestimate}) and the fact that $\frac{\leb^{2(d-1)}(F_{j,k,l})}{2^l} \lesssim 1$.
Since $S_{k,l}$ is summable this gives
\begin{align*}
\|T[f]\|_{L^q(L^r),K} &\lesssim \sum_{j=-\infty}^{\infty} 2^{j} \leb^{d}(E_j)^{\frac{1-\epsilon_1}{p}}
\\ &= \sum_{j=-\infty}^0 2^{j\epsilon_1}  \left(2^j  \leb^{d}(E_j)^{\frac{1}{p}}\right)^{1-\epsilon_1} + \sum_{j=1}^\infty 2^{-j\epsilon_1}  \left(2^j  \leb^{d}(E_j)^{\frac{1-\epsilon_1}{p(1+\epsilon_1)}}\right)^{1+\epsilon_1}
\\ &\lesssim \|f\|_{L^{p}}^{1-\epsilon_1} + \|f\|_{L^{\frac{p(1+\epsilon_1)}{1-\epsilon_1}}}^{1+\epsilon_1}
\\ &\lesssim \|f\|_{L^{\frac{p(1+\epsilon_1)}{1-\epsilon_1}}}^{1-\epsilon_1} + \|f\|_{L^{\frac{p(1+\epsilon_1)}{1-\epsilon_1}}}^{1+\epsilon_1}
\\ &=2\|f\|_{L^{p+\epsilon}},
\end{align*}
where we choose $\epsilon_1$ solving $p\frac{1+\epsilon_1}{1-\epsilon_1} = p+\epsilon$, and recall (\ref{normfone}) for the last equation.
\end{proof}

\section{The two-ends reduction} \label{twoendssection}

In order to obtain a favorable value of $p$ in Theorem \ref{maintheorem}, we will employ a version of the \emph{two-ends reduction} from \cite{wo2}, namely that it suffices to assume for $\lambda \leq \rho \leq \lambda^{\epsilon}$, $(\xi,x) \in F$, $z \in \rea^d$, that
\begin{equation} \label{twoends}
\int_0^1 \chi_{E \cap B(z,\rho)}(x + t(\xi + e_d))\ dt \lesssim \lambda \rho^{B\epsilon}
\end{equation}
where $B = \frac{1}{4p}$. 
To justify this reduction, we will use induction on scales via Claim \ref{iosclaim} below.  

\begin{claim} \label{iosclaim}
Suppose that $1 + \frac{d-1}{r} \geq \frac{d}{p}$ and that the inequality 
\begin{equation} \label{ioshypr}
\leb^d(E) \geq C_\epsilon \lambda^{p+\epsilon} 
\leb^{2(d-1)}(F)^{\frac{p}{q}} \ma(F)^{p\left(\frac{1}{r} - \frac{1}{q}\right)} \end{equation}
holds for 
\begin{equation} \label{ioshyprreq}
E \subset Q,\ \lambda \geq \lambda_0, \text{\ and\ } E,F,\lambda \text{\ satisfying\ } (\ref{EFlambdahyp}).
\end{equation}
Then, we also have for $\rho < 1$ 
\begin{equation} \label{iosconclr}
\leb^d(\tilde{E}) \geq \rho^{-\epsilon} C_\epsilon \tilde{\lambda}^{p+\epsilon} 
\leb^{2(d-1)}(\tilde{F})^{\frac{p}{q}} \ma(\tilde{F})^{p\left(\frac{1}{r} - \frac{1}{q}\right)} 
\end{equation}
when $\tilde{E}$ is contained in any sub-cube $\tilde{Q} \subset Q$ with side length $\rho$; 
$\tilde{\lambda} \geq \rho \lambda_0$; and $\tilde{E}, \tilde{F}, \tilde{\lambda}$ are as in (\ref{EFlambdahyp}) (i.e. $T[\chi_{\tilde{E}}] \geq \tilde{\lambda}$ on $\tilde{F}$).
\end{claim}

\begin{proof}
The proof is simply a change of variables, mapping $\tilde{Q}$ to $Q$.  
Let $z$ be the center of $\tilde{Q}$.
Then, letting $t_z e_d= \proj_{e_d}(z)$, $x_z = \proj_H(z)$, and $\tilde{t} = \frac{1}{\rho}\left(t - \left(t_z - \frac{\rho}{2}\right)\right)$, we see that
\begin{eqnarray*}
T[\chi_{\tilde{E}}](\xi,x)&=&\int_{t_z - \frac{\rho}{2}}^{t_z + \frac{\rho}{2}}
\chi_{\tilde{E}}(x + t(\xi + e_d))\ dt \\ &=&
\int_{t_z - \frac{\rho}{2}}^{t_z + \frac{\rho}{2}}
\chi_{\tilde{E} - z + \frac{\rho}{2}e_d}\left(x - x_z + t(\xi + e_d) - \left(t_z - \frac{\rho}{2}\right)e_d\right)\ dt
\\ &=&  \rho \int_{0}^{1}
\chi_{\tilde{E} - z + \frac{\rho}{2}e_d}\left(x - x_z  + \left(t_z - \frac{\rho}{2}\right)\xi + \rho \tilde{t}\left(\xi + e_d\right) \right)\ d\tilde{t} \\ &=&
\rho \int_{0}^{1}
\chi_{\frac{1}{\rho}\left(\tilde{E} - z + \frac{\rho}{2}e_d\right)}\left(\frac{1}{\rho}\left(x - x_z + \left(t_z - \frac{\rho}{2}\right) \xi \right) + \tilde{t}(\xi + e_d) \right)\ d\tilde{t} \\ &=&
\rho T\left[\chi_E \right]\left(\xi,\frac{1}{\rho}\left(x - x_z + \left(t_z - \frac{\rho}{2}\right) \xi \right)\right),
\end{eqnarray*}
where $E =\frac{1}{\rho}\left(\tilde{E} - z + \frac{\rho}{2}e_d\right)$.  

Let $F = \{(\xi,\frac{1}{\rho}\left(x - x_z + \left(t_z - \frac{\rho}{2}\right)\xi \right)) : (\xi,x) \in \tilde{F}\}$ and $\lambda = \frac{1}{\rho} \tilde{\lambda}$, so that $T[\chi_E] \geq \lambda$ on $F$.
Then
\begin{gather*}
\leb^d(E) = \rho^{-d} \leb^d(\tilde{E}),
\\ \leb^{2(d-1)}(F) = \rho^{-(d-1)} \leb^{2(d-1)}\tilde{F},
\intertext{and} \ma(F) = \rho^{-(d-1)} \ma(\tilde{F}).
\end{gather*}

By construction, $E$ is contained in $Q$ and,
 since $\tilde{\lambda} \geq \rho \lambda_0$, we have $\lambda \geq \lambda_0$.
Thus, we may apply (\ref{ioshypr}), obtaining 
\begin{eqnarray*}
\rho^{-d} \leb^{d}(\tilde{E}) &\geq& C_\epsilon \left(\rho^{-1}\tilde{\lambda}\right)^{p+\epsilon} 
\left(\rho^{-(d-1)}\leb^{2(d-1)}(\tilde{F})\right)^{\frac{p}{q}} (\rho^{-(d-1)}\ma(\tilde{F}))^{p\left(\frac{1}{r} - \frac{1}{q}\right)}
\\ &=& \rho^{-p(1 + \frac{d-1}{r})} \rho^{-\epsilon} C_{\epsilon} \tilde{\lambda}^{p+\epsilon} \leb^{2(d-1)}(\tilde{F})^{\frac{p}{q}}\ma(\tilde{F})^{p\left(\frac{1}{r} - \frac{1}{q}\right)}.
\end{eqnarray*}

Then, since $1 + \frac{d-1}{r} \geq \frac{d}{p}$ and $\rho < 1$, we have $\rho^{d -p(1 + \frac{d-1}{r})} > 1$ and
thus our conclusion (\ref{iosconclr}) holds. 
\end{proof}

We will use induction on $\lambda_0$ to show that proving (\ref{ioshypr})
for arbitrary $\lambda$ under the two-ends reduction is sufficient to 
prove (\ref{ioshypr}) for arbitrary $\lambda$ in general. 

Note that (\ref{twoends}) holds trivially for any fixed choice of $\lambda$ by a sufficiently large choice of the implicit constant. Thus, we may start our induction and assume that (\ref{ioshypr}) holds for an initial $\lambda_0$. 
Now, suppose that (\ref{ioshypr}) is known to hold with $\lambda_0 = \Lambda \leq 64^{-\frac{1}{\epsilon}}$ and that we want to prove it with $\lambda_0 = \frac{1}{2} \Lambda$. Let $\frac{1}{2}\Lambda\leq \lambda \leq \Lambda$, and let $E, F$ be as in (\ref{ioshyprreq}).
Let $\overline{F}$ be the subset of $F$ for which the two ends condition (\ref{twoends}) fails.
If 
\begin{equation} \label{twoendsfailureequation}
\leb^{2(d-1)}(\overline{F}) \leq \frac{1}{2} \leb^{2(d-1)}(F),
\end{equation}
then we may apply our knowledge of (\ref{ioshypr}) under the two-ends reduction
to the set $F \setminus \overline{F}$. Thus, without loss of generality 
we assume the inequality opposite to (\ref{twoendsfailureequation}).

By dyadic pigeonholing we may find $\rho_0 \in [\lambda, \lambda^\epsilon]$ so that (\ref{twoends}) fails with $\rho=\rho_0$ for $(\xi,x) \in \hat{F}$ and 
\[
\leb^{2(d-1)}(\hat{F}) \gtrsim \frac{\leb^{2(d-1)}(F)}{|\log(\lambda)|}.
\] 
We then tile $Q$ by a collection $\{Q_j\}$ of cubes with side length $\rho_0$, and let $\tilde{E}_j = E \cap \tilde{Q}_j$ where $\tilde{Q}_j$ is the cube with the same center as $Q_{j}$ and side length $4\rho_0$. For each $(\xi,x) \in \hat{F}$, (\ref{twoends}) fails for some $z_{\xi,x} \in Q_{j_{\xi,x}}$, giving 
\begin{equation} \label{usetefailure}
T[\chi_{\tilde{E}_{j_{\xi,x}}}](\xi,x) \geq \lambda \rho_0^{B\epsilon}.
\end{equation}
Henceforth, we take $\tilde{\lambda} = \lambda \rho_0^{B\epsilon}$ and let $\tilde{F}_j = \{(\xi,x): j = j_{\xi,x}\}$. 
Without loss of generality, we assume that the constant on the right hand side of (\ref{twoends}) is at least $1$ and thus we may ignore it in (\ref{usetefailure}).

For each $j$ we now have $T[\chi_{\tilde{E}_j}] \geq \tilde{\lambda}$ on $\tilde{F}_j$, and $\tilde{E_j}$ contained in a sub-cube of $Q$ with side length $4\rho_0$. Naturally, we will want to apply Claim \ref{iosclaim} in order to estimate the size
of the $\tilde{E}_j$. In order to apply the claim, it only remains to verify that $\tilde{\lambda} \geq 4\rho_0 \Lambda$.
In fact, assuming without loss of generality that $\epsilon \leq \frac{1}{2B}$,
we have 
\[
\tilde{\lambda} \geq  \lambda \rho_0^{B\epsilon} = 
\rho_0 \Lambda  \frac{\lambda}{\Lambda} \rho_0^{B\epsilon - 1} 
\geq \rho_0 \Lambda  \frac{1}{2} \rho_0^{B \epsilon - 1}
\geq \rho_0 \Lambda  \frac{1}{2} \rho_0^{-\frac{1}{2}}
\geq \rho_0 \Lambda  \frac{1}{2} \Lambda^{-\frac{1}{2}\epsilon}
\geq 4\rho_0 \Lambda.
\]
Thus, we may apply Claim \ref{iosclaim}. We then have by (\ref{iosconclr})
\begin{eqnarray*}
\leb^d(\tilde{E}_j) &\geq& (4\rho_0)^{-\epsilon} C_{\epsilon} \tilde{\lambda}^{p+\epsilon} 
\leb^{2(d-1)}(\tilde{F}_j)^{\frac{p}{q}} \ma(\tilde{F_j})^{p\left(\frac{1}{r} - \frac{1}{q}\right)} \\
&\geq&\frac{1}{4} \rho_0^{-\epsilon + (p+\epsilon)B\epsilon} C_{\epsilon} \lambda^{p+\epsilon}
\leb^{2(d-1)}(\tilde{F}_j)^{\frac{p}{q}} \ma(\tilde{F_j})^{p\left(\frac{1}{r} - \frac{1}{q}\right)}.
\end{eqnarray*}

Recalling that $q \leq r$, we observe that, since $\tilde{F}_j \subset F$, we have $\ma(\tilde{F}_j)^{p(\frac{1}{r}-\frac{1}{q})} \geq \ma(F)^{p(\frac{1}{r}-\frac{1}{q})}$. Since the $\tilde{E}_j$ are finitely overlapping and $q \geq p$, we then have
\begin{eqnarray} \label{implicitconstantov}
\leb^d(E) &\gtrsim& \sum_{j}\leb^d(\tilde{E}_j) 
\\ \notag &\geq& \frac{1}{4}\rho_0^{-\epsilon + (p+\epsilon)B\epsilon}  C_{\epsilon} \lambda^{p+\epsilon} 
\sum_j \leb^{2(d-1)}(\tilde{F}_j)^{\frac{p}{q}}\ma(\tilde{F}_j)^{p\left(\frac{1}{r} - \frac{1}{q}\right)}
\\ \nonumber &\geq& \frac{1}{4}\rho_0^{-\epsilon + (p+\epsilon)B\epsilon} C_{\epsilon} \lambda^{p+\epsilon} 
\ma(F)^{p\left(\frac{1}{r} - \frac{1}{q}\right)}
\sum_j \leb^{2(d-1)}(\tilde{F}_j)^{\frac{p}{q}}
\\ \nonumber &\geq& \frac{1}{4}\rho_0^{-\epsilon + (p+\epsilon)B\epsilon} C_{\epsilon} \lambda^{p+\epsilon}  \ma(F)^{p\left(\frac{1}{r} - \frac{1}{q}\right)} \leb^{2(d-1)}(\hat{F})^{\frac{p}{q}}
\\ \nonumber &\gtrsim& \rho_0^{-\epsilon + (p+\epsilon)B\epsilon} |\log(\lambda)|^{-\frac{p}{q}} C_{\epsilon} \lambda^{p+\epsilon} 
\ma(F)^{p\left(\frac{1}{r} - \frac{1}{q}\right)} \leb^{2(d-1)}(F)^{\frac{p}{q}}.
\end{eqnarray}
Recalling that $B = \frac{1}{4p}$ and assuming that $\epsilon < p$, we have 
\begin{equation} \label{lastioseq}
\rho_0^{-\epsilon + (p+\epsilon)B\epsilon} |\log(\lambda)|^{-\frac{p}{q}}
\geq
\lambda^{-\frac{1}{2}\epsilon^2} |\log(\lambda)|^{-\frac{p}{q}}.
\end{equation}
It only remains to verify that the right side of (\ref{lastioseq}) is large enough to overcome 
the implicit constant in (\ref{implicitconstantov}). Noting that this 
constant is independent of the constant in (\ref{twoends}), we see that this 
may be accomplished by a sufficiently small initial choice of $\lambda_0$.

\section{The main estimate} \label{continslice}

For $t \in \rea$, let $H_t$ denote the plane $H + t e_d$. Given a line $g$ which intersects $H$ in exactly one point, we define
\[
\pi_t(g) = H_t \cap g.
\] 
For collections $G$ of such lines, we are interested in lower bounds for the size of $\pi_t(G)$ in terms of the size of $G$. Suppose $G$ is finite. Then, from the fact that 
a line is determined by two points, we see that for every $t_1 \neq t_2$,
\[
\#G^{\frac{1}{2}} \leq \sup_{t = t_1,t_2} \#\pi_t (G)
\]
where $\#$ denotes cardinality. In \cite{kt}, Katz and Tao showed that if $t_0,t_\infty,t_{1},t_{1'},t_{2},t_{2'}$ satisfy a certain 
algebraic condition, and if the lines in $G$ point in distinct directions, then
\begin{equation} \label{katztaofinite}
\#G^{\frac{4}{7}} \lesssim \sup_{t = t_0,t_\infty,t_{1},t_{1'},t_{2},t_{2'}} \#\pi_t (G).
\end{equation}
In order to give a bound of type (\ref{kakeyaorder}), one must consider the more general case where the lines in $G$ do not point in distinct directions. Suppose that at most $M$ lines in $G$ point in each direction. Then by taking a maximal direction separated subset of $G$, the estimate  
\begin{equation} 
\#G^{\frac{4}{7}} \lesssim M^{\frac{4}{7}} \sup_{t = t_0,t_\infty,t_{1},t_{1'},t_{2},t_{2'}} \#\pi_t (G)
\end{equation}
follows trivially from (\ref{katztaofinite}). However, following the proof of (\ref{katztaofinite}) with the quantity $M$ in mind, one actually obtains
\begin{equation} \label{correctMexponent}
\#G^{\frac{4}{7}} \lesssim M^{\frac{1}{7}} \sup_{t = t_0,t_\infty,t_{1},t_{1'},t_{2},t_{2'}} \#\pi_t (G)
\end{equation}
with no additional arguments required. This is, in fact, the sharp power of $M$
for the given power of $G$, as one may verify by letting $G_n$ be the set of lines determined by the pairs of points $(i e_1, i e_1 + j e_1 + e_d)$ where $1 \leq i,j \leq n$, letting $t_0,t_\infty,t_1,t_{1'},t_2,t_{2'}$ be any fixed rational numbers, and considering $n$ large.   

The maximal operator bound in \cite{kt} was proven using a $\delta$-discrezation argument and a refinement of (\ref{katztaofinite}) which took into account possible $\delta$-uncertainties.
It is possible to adapt (\ref{correctMexponent}) for use with a discretization argument and use this estimate to prove Theorem \ref{maintheorem}. Instead, we will prove the analog of (\ref{correctMexponent}) for Lebesgue-measurable sets of lines and avoid discretization entirely.  

Let $G$ be a set of lines in $\rea^d$, each of which intersect the plane $H$ in 
exactly one point. Considering the coordinates for the $x$-ray transform $T$, 
we parametrize $G$ by the subset 
\[
G_X = \{(\xi,x) : \text{\ there\ exists\ } g \in G \text{\ with\ } x,(x+\xi+e_d)\in g\}
\]
of $H \times H$.
For $t_1 \neq t_2$, we may also parametrize $G$ by the subset 
\[
G_{t_1,t_2} = \{(\pi_{t_1}(g),\pi_{t_2}(g)) : g \in G\}
\]
of $H \times H$.
Using the ``line property'' 
\begin{equation} \label{lineproperty}
(x,y) \in G_{t_1,t_2} \Leftrightarrow \left(
\frac{t_2 - t_3}{t_2 - t_1}x + \frac{t_3 - t_1}{t_2 - t_1}y,
 \frac{t_2 - t_4}{t_2 - t_1}x + \frac{t_4 - t_1}{t_2 - t_1}y 
\right) \in G_{t_3,t_4},
\end{equation}
we may change variables to give
\begin{equation} \label{measG}
|G| := \leb^{2(d-1)}(G_X) = \leb^{2(d-1)}(G_{0,1}) = |t_1 - t_2|^{-(d-1)} \leb^{2(d-1)}(G_{t_1,t_2}). 
\end{equation}
Henceforth, we will use the abbreviation
\[
D_{t,t'} := |t - t'|^{-(d-1)}.
\]
Finally, define 
\[
\ma(G) := \ma(G_X) = \sup_{\xi \in H} \leb^{d-1}(\{x : (\xi,x) \in G_X\}).
\]
In practice, each line in $G$ will intersect the cube $Q$, and $G_X \subset 
B(0,C) \times \rea^{d-1}.$ Thus $\ma(G),|G|,\leb^{d-1}(\pi_t(G)) < \infty$. We will 
always assume that this inequality holds below. To put the following proposition
in context, we point out that, as we will see later, the quantity $\frac{k}{\alpha - \beta}$ in an estimate of the form (\ref{iterativemainestimatehyp}) corresponds to the quantity $\frac{r}{p}$ in an estimate of the form (\ref{kakeyaorder}).

\begin{proposition} \label{iterativemainestimate}
Suppose that for $t_1,\ldots,t_k \in \rea$, we have 
\begin{equation} \label{iterativemainestimatehyp}
|G|^{\alpha} \lesssim C_{t_1,\ldots,t_k} \ma(G)^{\beta}  \prod_{i = 1, \ldots, k} \leb^{d-1}(\pi_{t_i}(G))
\end{equation}
for each set of lines $G$, where $0 \leq \beta$ and $\alpha \leq k$.

Then, for any $t_0, t_\infty, t_{1'}, \ldots, t_{k'} \in \rea$
satisfying $t_0 \neq t_\infty$;
\[
t_i \neq t_0,\ \ t_{i'} \neq t_0,\ \ t_i \neq t_\infty, \text{\ for\ }i = 1, \ldots, k;
\]
and satisfying the requirement that 
\begin{equation} \label{iteratetildesdef}
s = s_{t_0,t_\infty}(t_i,t_{i'}) := (t_\infty - t_0) \frac{t_i - t_0}{(t_i - t_\infty)(t_{i'} - t_0)}
\end{equation}
is independent of $i$ for $i = 1, \ldots, k$,
we have
\begin{equation} \label{iterativemainestimateconclusion}
\begin{split}
|G|^{2k} &\lesssim \left(C_{t_1, \ldots, t_k} D_{t_0,t_\infty}^{k-\alpha} \prod_{i = 1, \ldots k} D_{t_{i'},t_0}\right) \ma(G)^{\beta + k - \alpha}  
\\ &\quad \cdot \leb^{d-1}(\pi_{t_\infty}(G))^{k-\alpha} \leb^{d-1}(\pi_{t_0}(G))^{k} \prod_{i = 1, \ldots k} \leb^{d-1}(\pi_{t_i}(G)) \leb^{d-1}(\pi_{t_{i'}}(G)).  
\end{split}
\end{equation}
\end{proposition}

For $t_1 \neq t_2$, the trivial estimate
\begin{equation} \label{rpequal2}
|G| \leq D_{t_1,t_2} \leb^{d-1}(\pi_{t_1}(G)) \leb^{d-1}(\pi_{t_2}(G))
\end{equation}
follows directly from (\ref{measG}). Thus, applying Proposition \ref{iterativemainestimate} with $k = 2, \alpha = 1, \beta = 0,$ and $C_{t_1,t_2} = D_{t_1,t_2}$, we obtain

\begin{corollary} \label{sixsliceestimate}
Let $t_1, t_{1'}, t_{2}, t_{2'}, t_0, t_{\infty} \in \rea$ satisfy
\[
t_1 \neq t_2;\ t_0 \neq t_{\infty}, t_1, t_{1'}, t_2, {t_2'};\ t_{\infty} \neq t_1,t_2;   
\]
and the requirement that $s$ is independent of $i$ in (\ref{iteratetildesdef}).
Then
\begin{align} \label{sixslicestatement}
|G|^{4} &\lesssim \left(D_{t_1,t_2} D_{t_0,t_\infty} D_{t_{1'},t_0}, D_{t_2',t_0} \right)  \ma(G)
\\ \notag &\ \ \ \ \ \ \ \cdot \sup_{t = t_{\infty}, t_0, t_1,t_{1'},t_2,t_{2'}} \leb^{d-1}(\pi_{t}(G))^{7}.
\end{align}
\end{corollary}

In Section \ref{pigeonholesection}, we will use a uniformization argument to 
obtain Theorem \ref{maintheorem} from Corollary \ref{sixsliceestimate}.
In Section \ref{pigeonholesection2}, we will consider further iterations of Proposition \ref{iterativemainestimate}.

\begin{proof}[Proof of Proposition \ref{iterativemainestimate}]
We will begin by defining the set 
\[
V = \{(g_1,g_2) \in G \times G : \pi_{t_0}(g_1) = \pi_{t_0}(g_2)\}.
\]
For any $w = (g_1,g_2) \in V$, let $\gamma_i(w) = g_i$ when $i = 1,2$.
Consider the function on $V$,
\[
\nu = s \pi_{t_\infty}(\gamma_1) + \pi_{t_\infty}(\gamma_2) - \pi_{t_0}(\gamma_2).
\]
Our purpose in defining $\nu$ is to obtain the equivalence classes in subsets $W$ of $V$ determined by the fibers of $\nu$: 
\[
W_{\nu_{0}} := W \cap \nu^{-1}(\nu_0) \text{\ for\ } \nu_0 \in H.
\]

Let $G_{\nu_0} = \gamma_1(V_{\nu_0}).$
Then we may calculate an upper bound for $|G_{\nu_0}|$ in terms 
of $\ma(G)$ and $\leb^{d-1}(\pi_{t_{\infty}}(G)).$
First note that
\begin{align} \label{gnuubone}
|G_{\nu_0}| &= D_{t_0,t_\infty} \int_{\rea^{d-1}} \int_{\rea^{d-1}} \chi_{(G_{\nu_{0}})_{t_0,t_{\infty}}}(x,y)\ dx\ dy
\\ \notag &= D_{t_0,t_\infty} \int_{\rea^{d-1}} \int_{\rea^{d-1}} \chi_{G_{t_0,t_{\infty}}}(x,y) \chi_{G_{t_0,t_\infty}}(x,\nu_{0} - sy + x)\ dx\ dy.
\end{align}
Fix $y$ and let $\xi_y = \frac{\nu_0 - s y}{t_\infty - t_0}$. Observe that
\begin{align} \label{gnuubtwo}
\int_{\rea^{d-1}} \chi_{G_{t_0,t_\infty}}(x,x + \nu_{0} - sy)\ dx &= 
\int_{\rea^{d-1}} \chi_{G_{t_0,t_\infty}}(x,x + (t_\infty - t_0)\xi_y)\ dx
\\ \notag &=\int_{\rea^{d-1}} \chi_{G_{0,1}}(x - t_0 \xi_y, x  - t_0 \xi_y + \xi_y)\ dx
\\ \notag &= \leb^{d-1}(\{x : (\xi_y,x) \in G_X\})
\\ \notag &\leq \ma(G),
\end{align}
where the second equation follows from the line property (\ref{lineproperty}), 
and the third equation follows from the translation invariance of $\leb^{d-1}$.
Combining (\ref{gnuubone}) and (\ref{gnuubtwo}), we have
\begin{equation} \label{gnuub}
|G_{\nu_0}| \leq D_{t_0,t_\infty} \leb^{d-1}(\pi_{t_\infty}(G)) \ma(G).
\end{equation}

The remainder of our argument will consist of using (\ref{iterativemainestimatehyp}) to obtain a 
lower bound for $|G_{\nu_0}|$ when $\nu_0$ is chosen favorably.
For any $t,t' \neq t_0$ and subset $W$ of $V$, we may parametrize $W$ by the subset
\[
W_{t_0,t,t'} = \{(\pi_{t_0}(g_1),\pi_{t}(g_1), \pi_{t'}(g_2)) : (g_1,g_2) \in W \}
\]
of $H^3$.
Using our line-property (\ref{lineproperty}), we note that, as was the case
with $|G|$,  
\[
|W| := D_{t_0,t} D_{t_0,t'} \leb^{3(d-1)}(W_{t_0,t,t'})
\]
is independent of $t,t'$.

For any $t,t' \neq t_0$, $W \subset V$, we may consider the subset of $W$ 
which is ``popular'' with respect to the double projection $(\pi_{t}(\gamma_1),\pi_{t'}(\gamma_2))$
\begin{equation}
\begin{split}
W^{\langle \pi_{t \otimes t'} \rangle}_{t_0,t,t'} := 
\left\{(x,y,z) : D_{t_{0},t} D_{t_0,t'} \int_{\rea^{d-1}} \chi_{W_{t_0,t,t'}}(x',y,z) \ dx' \right.
\\ \left. \geq \frac{|W|}{2\leb^{d-1}(\pi_{t}(G))\leb^{d-1}(\pi_{t'}(G))}\vphantom{\int}\right\}.
\end{split}
\end{equation}
After estimating $\left|W \setminus W^{\langle \pi_{t \otimes t'} \rangle}\right|$, one observes that
\begin{equation}
|W^{\langle \pi_{t \otimes t'} \rangle}| \geq \frac{1}{2} |W|.
\end{equation}
Thus, abbreviating
\[
V' := \left(\left(\left(V^{\langle \pi_{t_1 \otimes t_1'} \rangle}\right)^{\langle \pi_{t_2 \otimes t_2'} \rangle}\right)...\right)^{\langle \pi_{t_k \otimes t_k'} \rangle},
\]
we have
\begin{equation} \label{refinementequation}
|V'| \gtrsim |V|.
\end{equation}

Given any subset $W$ of $V$,  
\begin{align} \label{integratenu}
|W| &= D_{t_{0},t_{\infty}} D_{t_0,t_\infty} \int_{\rea^{d-1}} \int_{\rea^{d-1}}\int_{\rea^{d-1}} \chi_{W_{t_0,t_\infty,t_\infty}}(x,y,z)\ dx\ dy\ dz
\\ \notag &= D_{t_{0},t_{\infty}} D_{t_0,t_\infty} \int_{\rea^{d-1}} \int_{\rea^{d-1}}\int_{\rea^{d-1}} \chi_{W_{t_0,t_\infty,t_\infty}}(x,y,\nu' - s y + x)\ dx\ dy\ d\nu'
\\ \notag &= D_{t_{0},t_{\infty}} D_{t_0,t_\infty} \int_{\rea^{d-1}} \int_{\rea^{d-1}}\int_{\rea^{d-1}} \chi_{(\gamma_1(W_{\nu'}))_{t_0,t_\infty}}(x,y)\ dx\ dy\ d\nu'.
\end{align}
Substituting (\ref{integratenu}) with $W = V'$ and $W = V$ into the left and right hand sides respectively of (\ref{refinementequation}), we observe that
\begin{equation} \label{gnulbone}
|G'_{\nu_0}| \gtrsim |G_{\nu_0}|
\end{equation}
for some $\nu_{0} \in \nu(V)$, where we define 
$G'_{\nu_0} = \gamma_1\left((V')_{\nu_0}\right)$.

Applying our hypothesis (\ref{iterativemainestimatehyp}) to the set of lines $G'_{\nu_0}$, we obtain
\begin{align} \label{gnulbtwo}
|G'_{\nu_0}|^{\alpha} &\lesssim C_{t_1,\ldots,t_k} \ma(G'_{\nu_0})^{\beta} 
\prod_{i = 1, \ldots, k}\leb^{d-1}(\pi_{t_{i}}(G'_{\nu_0}))
\\ \notag &\leq C_{t_1,\ldots,t_k} \ma(G)^{\beta} 
\prod_{i = 1, \ldots, k}\leb^{d-1}(\pi_{t_{i}}(G'_{\nu_0})),
\end{align}
where the second inequality follows from the fact that $G'_{\nu_0} \subset G$ and the condition that $\beta \geq 0$.

Suppose $y \in \pi_{t_{i}}(G'_{\nu_0}).$ Then 
$y = \pi_{t_{i}}(g_1)$ where $(g_1,g_2) \in (V')_{\nu_0}$. Letting 
$z = \pi_{t_{i'}}(g_2)$, we have by definition 
of $V'$
\begin{equation} \label{usingdeterminedby1}
\frac{|V|}{\leb^{d-1}(\pi_{t_i}(G))\leb^{d-1}(\pi_{t_{i'}}(G))}
\lesssim D_{t_i,t_0} D_{t_{i'},t_0} \int_{\rea^{d-1}}
\chi_{V_{t_{0},t_i,t_{i'}}}(x,y,z)\ dx.
\end{equation}
We will now need to use our definition of $s$,
 rewriting
\begin{align} \label{determinesnu}
\nu &= s \left(\pi_{t_\infty}(\gamma_1) - \pi_{t_0}(\gamma_1)\right)
+ \left(\pi_{t_\infty}(\gamma_2) - \pi_{t_0}(\gamma_2)\right)
+ s \pi_{t_0}(\gamma_2)
\\ \notag &= s \frac{t_\infty - t_0}{t_i - t_0}\left(\pi_{t_i}(\gamma_1) - \pi_{t_0}(\gamma_1)\right)
 + \frac{t_\infty - t_0}{t_{i'} - t_0}\left(\pi_{t_{i'}}(\gamma_2) - \pi_{t_0}(\gamma_2)\right) 
+ s \pi_{t_0}(\gamma_2)
\\ \notag &= s \frac{t_\infty - t_0}{t_i - t_0} \pi_{t_i}(\gamma_1)  
+ \frac{t_\infty - t_0}{t_{i'} - t_0} \pi_{t_{i'}}(\gamma_2) 
+ \pi_{t_0}(\gamma_2)\left(\frac{t_{0} - t_\infty}{t_{i'} - t_0} + s \frac{t_i - t_\infty}{t_i - t_0}\right)
\\ \notag &= s \frac{t_\infty - t_0}{t_i - t_0} \pi_{t_i}(\gamma_1)  + \frac{t_\infty - t_0}{t_{i'} - t_0} \pi_{t_{i'}}(\gamma_2),
\end{align}
where $\pi_{t_0}(\gamma_j)$ is independent of $j$ by definition of $V$,
and where the last equation follows from (\ref{iteratetildesdef}). The point is that membership in $V_{\nu_0}$ is determined by the double projection $(\pi_{t_i}(\gamma_1),\pi_{t_i'}(\gamma_2))$. Hence
\begin{align} \label{usingdeterminedby2}
\int_{\rea^{d-1}}
\chi_{V_{t_{0},t_i,t_{i'}}}(x,y,z)\ dx 
&=  \int_{\rea^{d-1}}
\chi_{(V_{\nu_0})_{t_{0},t_i,t_{i'}}}(x,y,z)\ dx
\\ \notag &= \int_{\rea^{d-1}}
\chi_{(G_{\nu_0})_{t_{0},t_i}}(x,y)\ dx.
\end{align}

Combining (\ref{usingdeterminedby1}) and (\ref{usingdeterminedby2}), we obtain
\begin{align} \label{gnulbthree}
\leb^{d-1}(\pi_{t_i}(G'_{\nu_0})) 
&=  
\int_{\rea^{d-1}} \chi_{\pi_{t_i}(G'_{\nu_0})}(y)
\ dy
\\ \notag &\lesssim 
D_{t_i,t_0} D_{t_{i'},t_0} \frac{\leb^{d-1}(\pi_{t_i}(G))\leb^{d-1}(\pi_{t_{i'}}(G))}{|V|} 
\\ \notag & \ \ \ \ \ \ \ \ \ \cdot \int_{\rea^{d-1}} \int_{\rea^{d-1}} 
\chi_{(G_{\nu_0})_{t_0,t_i}}(x,y)\ dx\ dy 
\\ \notag &=
|G_{\nu_0}| D_{t_{i'},t_0} \frac{\leb^{d-1}(\pi_{t_i}(G))\leb^{d-1}(\pi_{t_{i'}}(G))}{|V|}.
\end{align}
Combining (\ref{gnulbone}), (\ref{gnulbtwo}), and (\ref{gnulbthree}),
we have 
\begin{align} \label{gnulbfour}
|V|^{k} &\lesssim C_{t_1,\ldots,t_k,\tau} \ma(G)^{\beta} |G_{\nu_0}|^{k-\alpha} 
\\ \notag &\ \ \ \ \ \ \ \ \cdot \prod_{i=1,\ldots,k} D_{t_0,t_i'} \leb^{d-1}(\pi_{t_i}(G)) \leb^{d-1}(\pi_{t_{i'}}(G)).
\end{align}

From the Cauchy-Schwarz inequality, we obtain a lower bound for $|V|$ in terms
of $|G|$:
\begin{align*}
\left(\frac{|G|}{D_{t_0,t}}\right)^2 &\leq \leb^{d-1}(\pi_{t_0}(G)) \int_{\rea^{d-1}} \int_{\rea^{d-1}} \chi_{G_{t_0,t}}(x,y)\ dy \int_{\rea^{d-1}} \chi_{G_{t_0,t}}(x,z)\ dz\ dx
\\ &= \leb^{d-1}(\pi_{t_0}(G)) \int_{\rea^{d-1}} \int_{\rea^{d-1}} \int_{\rea^{d-1}} \chi_{V_{t_0,t,t}}(x,y,z)\ dx\ dy\ dz
\\ &= \leb^{d-1}(\pi_{t_0}(G)) \frac{|V|}{D_{t_0,t}^2},
\end{align*}
which simplifies to
\begin{equation} \label{cauchyschwarzv}
|V| \geq \frac{|G|^2}{\leb^{d-1}(\pi_{t_0}(G))}.
\end{equation}

Combining (\ref{cauchyschwarzv}), (\ref{gnulbfour}), and (\ref{gnuub}), we finally obtain
(\ref{iterativemainestimateconclusion}).
\end{proof}

\section{Uniformization} \label{pigeonholesection}
We now want to find six suitably ``average'' slices of $E$ to which we will 
apply Corollary \ref{sixsliceestimate}, allowing us to prove (\ref{ioshypr})
under the two-ends reduction. Our argument below is a continuous version of the uniformization argument in \cite{kt}.
Let $E,F,\lambda$ be as in (\ref{ioshyprreq}), and define
\begin{eqnarray*}
\mass(E,F) &=& \int_{\rea^{d-1}}\int_{\rea^{d-1}}\chi_{F}(\xi,x)T[\chi_{E}](\xi,x)\ dx\ d\xi
\\ &=&\int_{\rea^{d-1}}\int_{\rea^{d-1}}\int_0^1 \chi_F(\xi,x) \chi_{E}(x + t(\xi + e_d))\ dt\ dx\ d\xi 
\\ &=& \int_{\rea^{d}}\chi_{E}(z) \int_{\rea^{d-1}}\chi_{F}(\xi, x_z - t_z \xi)\ d\xi\ dz,
\end{eqnarray*}
where $x_z = \proj_{H}(z)$ and $t_z e_d = \proj_{e_d}(z).$
By definition of $F$, $\mass(E,F) \geq \lambda \leb^{2(d-1)}(F)$.

We will abbreviate $\gamma_{d} = \proj_{e_d}$. Let
\[
S_0 = \{t \in [0,1] : \leb^{d-1}(E \cap \gamma_d^{-1}(t)) \ll \lambda \leb^{2(d-1)}(F) \}
\]
and for $k > 0$ let  
\[
S_k = \{t \in [0,1]: (\lambda \leb^{2(d-1)}(F))^{1-(k-1)\epsilon} \lesssim \leb^{d-1}(E \cap \gamma_d^{-1}(t)) \ll (\lambda \leb^{2(d-1)}(F))^{1-k\epsilon}\}.
\]
Recalling that $F \subset B(0,C) \times \rea^{d-1}$, we note that
\[
\int_{\rea^{d-1}}\chi_{F}(\xi, \proj_{H}(z) - \proj_{e_d}(z)\xi)\ d\xi \lesssim 1
\]
for every $z$. 
Hence, defining $E_k = E \cap \gamma_d^{-1}(S_k)$, we have $\mass(E_0,F) \ll \lambda \leb^{2(d-1)}(F)$.
Thus, since $E = \bigcup_{k \lesssim \frac{1}{\epsilon}} E_k$, an appropriate choice of implicit constants gives 
\[
\mass(E_{k_0},F) \gtrsim \frac{1}{\epsilon} \lambda \leb^{2(d-1)}(F)
\gtrsim \lambda \leb^{2(d-1)}(F)
\]
for some $k_0 > 0.$ 
Let $E' = E_{k_0}$, $S = S_{k_0}$, and 
\[
F' = \{(\xi,x) \in F : T[\chi_{E'}](\xi,x) \gtrsim \lambda \}.
\]
Considering $\mass(E,F \setminus F')$, we note that
\[
\mass(E',F') \gtrsim \lambda \leb^{2(d-1)}(F).
\] 

We now proceed to find a point $(t_0,t_{\infty},t_1,t_{1'},t_2,t_{2'}) \in S^{6}$
with which we may apply 
Corollary \ref{sixsliceestimate} to our advantage. Due to the factors of $D_{t_i,t_j}$ in 
(\ref{sixslicestatement}), we would like to keep $|t_i - t_j|$ suitably large; 
this is facilitated by the two-ends reduction.  

For every $(\xi,x) \in F'$, let 
\[
S_{\xi,x} = \{t' \in S : x + t'(\xi + e_d) \in E'\}.
\]
Then, by definition of $F'$, 
\[
\mu_{\xi,x} := \leb^{1}(S_{\xi,x}) \gtrsim \lambda.
\]
By (\ref{twoends}) we have, for every $t \in \rea,$
\[
\leb^{1}(\{t' \in S_{\xi,x} : |t' - t| < \lambda^{\epsilon}\}) \lesssim \lambda \lambda^{B\epsilon^2} \ll \lambda
\]
where we assume, without loss of generality, that $\lambda$ is sufficiently small to obtain the rightmost inequality. Thus,
\[
\leb^{1}(\{t' \in S_{\xi,x} : |t' - t| \geq \lambda^{\epsilon}\}) \gtrsim \mu_{\xi,x}
\] 
and so letting  
\[
P(\xi,x) = \{(t_0,t_\infty) \in (S_{\xi,x})^2 : |t_0 - t_\infty| \geq \lambda^{\epsilon}\},
\]
we have $\leb^{2}(P(\xi,x)) \gtrsim \mu_{\xi,x}^{2}$. 

For each $(t_0,t_\infty) \in P(\xi,x)$ let 
\[
Q_{t_0,t_\infty}(\xi,x) = \{(t_1,t_{1'}) \in (S_{\xi,x})^2 : \text{\ for all \ } i \neq j \in \{0,\infty,1,1'\},\ |t_i - t_j| \geq \lambda^{\epsilon}
\},
\]
and note that 
\begin{equation} \label{measQequation}
\leb^{2}(Q_{t_0,t_\infty}(\xi,x)) \gtrsim \mu_{\xi,x}^{2}
\end{equation}
for every $(t_0,t_\infty) \in P(\xi,x)$.

We recall the definition of $s$ from (\ref{iteratetildesdef}), 
\[
s_{t_0,t_\infty}(t_i,t_{i'}) := \frac{(t_\infty - t_0)(t_i - t_0)}{(t_i - t_\infty)(t_{i'} - t_{0})}.
\]
In order to satisfy the condition in Corollary \ref{sixsliceestimate} that
$s$ is independent of $i$, we consider the set
\[
R_{t_0,t_\infty}(\xi,x) = \{(t_1,t_{1'},t_2,t_{2'}) \in Q_{t_0,t_\infty}(\xi,x)^2 : s_{t_0,t_\infty}(t_1,t_{1'}) = s_{t_0,t_\infty}(t_2,t_{2'})\}.
\]
Below we abbreviate $s_{t_0,t_\infty}$ by $s$, $Q_{t_0,t_\infty}(\xi,x)$ by $Q$,
$R_{t_0,t_\infty}(\xi,x)$ by $R$, and $\mu_{\xi,x}$ by $\mu$. Also we use $\cdot \gtrapprox \cdot$ to denote $\cdot \gtrsim \lambda^{C\epsilon} \cdot$, 
and we similarly use $\lessapprox$ and $\approxeq$.
We observe that $s(t_1,\cdot)$ is a diffeomorphism on an open set containing 
the support of $\chi_{Q(t_1,\cdot)}$, and so we may change variables to obtain
\begin{align} \label{changevarQequation}
\leb^{2}(Q) &= \int_{\rea} \int_{\rea} \chi_{Q}(t_{1},t_{1'})\ dt_1\ dt_{1'} 
\\ \notag &= \int_{\rea} \int_{\rea} \chi_{Q}\left(t_1,u_{t_0,t_\infty}(t_1,s')\right)\left|\frac{(u_{t_0,t_\infty}(t_1,s') - t_0)^2(t_1 - t_\infty)}{(t_1 - t_0)(t_\infty - t_0)} \right|\ dt_1\ ds'
\\ \notag &\approxeq \int_{\rea} \int_{\rea} \chi_{Q}\left(t_1,u_{t_0,t_\infty}(t_1,s')\right)\ dt_1\ ds'
\end{align}
where we define
\[
u_{t_0,t_\infty}(t_1,s) = t_0 + \frac{(t_\infty - t_0)(t_1 - t_0)}{s(t_1 - t_\infty)}
\]
and where the $\approxeq$ follows from the fact that the Jacobian is $\approxeq 1$ on $Q$.
We apply Cauchy-Schwarz and change variables again to see that 
\begin{align} \label{Qcsequation}
\lefteqn{\int_{\rea} \int_{\rea} \chi_{Q}\left(t_1,u_{t_0,t_\infty}(t_1,s')\right)\ dt_1\ ds'} \\ \notag
&\leq \leb^{1}(s(Q))^{\frac{1}{2}}
\left(\int_{\rea} \int_{\rea} \int_{\rea} 
\chi_{Q}\left(t_1,u_{t_0,t_\infty}(t_1,s')\right)
\chi_{Q}\left(t_2,u_{t_0,t_\infty}(t_2,s')\right)
\ dt_1\ dt_2\ ds'\right)^{\frac{1}{2}}
\\ \notag &\lessapprox
\left(\int_{\rea} \int_{\rea} \int_{\rea} 
\chi_{Q}\left(t_1,t_{1'}\right)
\chi_{Q}\left(t_2,u_{t_0,t_\infty}(t_2,s(t_1,t_{1'}))\right)
\ dt_1\ dt_{1'}\ dt_2\right)^{\frac{1}{2}}.
\end{align}
Abbreviating 
\[
w_{t_0,t_\infty}(t_1,t_1',t_2) = u_{t_0,t_\infty}(t_2,s(t_1,t_{1'}))
\]
we have by construction
\[
s(t_1,t_{1'}) = s(t_2,w_{t_0,t_\infty}(t_1,t_{1'},t_2)) 
\]
and hence
\[
\chi_R(t_1,t_{1'},t_2,w_{t_0,t_\infty}(t_1,t_{1'},t_{2})) = \chi_{Q}(t_1,t_{1'}) \chi_Q(t_2,w_{t_0,t_\infty}(t_1,t_{1'},t_2)).
\]
Thus, we combine (\ref{measQequation}), (\ref{changevarQequation}), and (\ref{Qcsequation}), to obtain 
\begin{eqnarray*}
|R| &:=& \int_{\rea^3} \chi_R(t_1,t_{1'},t_2, w_{t_0,t_\infty}(t_1,t_{1'},t_2))\ dt_1\ dt_{1'}\ dt_{2}
\\ &\gtrapprox& \mu^4.
\end{eqnarray*}

We have control over all of the $|t_i - t_j|$ relevant to Corollary \ref{sixsliceestimate}, 
except for $|t_1 - t_2|$.
However, since
\[
\begin{split}
\int_{\rea^3} \chi_R(t_1,t_{1'},t_2, w_{t_0,t_\infty}(t_1,t_{1'},t_2))\chi_{[0,r]}(|t_1 - t_2|)\ dt_1\ dt_{1'}\ dt_2
\\ \lesssim r \leb^{1}(\proj_{t_1}(R))\leb^{1}(\proj_{t_{1'}}(R)), 
\end{split}
\]
and 
$\leb^{1}(\proj_{t_i}(R)) \leq \mu$ for each $i$,
we have $|R'| \gtrapprox \mu^{4}$ where 
\begin{equation} \label{R'def}
R' = \{(t_1,t_{1'},t_2,t_{2'}) \in R : |t_1 - t_2| \gtrapprox \mu^2\}.
\end{equation}

Let
\begin{eqnarray*}
\lefteqn{X = \{(\xi,x,t_0,t_{\infty},t_1,t_{1'},t_2) : (\xi,x) \in F', (t_0,t_\infty) \in P(\xi,x), \text{\ and\ }} &&
\\ &&
\ \ \ (t_1,t_{1'},t_2,w_{t_0,t_\infty}(t_1,t_{1'},t_2)) \in R'_{t_0,t_\infty}(\xi,x)   
\}.
\end{eqnarray*}
Integrating everything out, we see that 
\[
\leb^{2(d-1)+5}(X) \gtrapprox \mass(E',F') \left(\inf_{(\xi,x) \in F'}\mu_{\xi,x}\right)^{5} \gtrsim \leb^{2(d-1)}(F) \lambda^{6}.
\]

For each $k \geq 0$ let 
\[
X_k = \{(\xi,x,t_0,t_\infty,t_1,t_1',t_2) :  \lambda^{C\epsilon} \lambda^{2-k\epsilon} \lesssim |t_1 - t_2| \ll \lambda^{C\epsilon} \lambda^{2-(k+1)\epsilon}\}.
\]
Then, recalling that each $\mu_{\xi,x} \gtrsim \lambda$ we have by definition of $R'_{\xi,x}$
\[
X = \bigcup_{k \leq \frac{2}{\epsilon}} X_{k}, 
\]
and hence for some $X' := X_{k_0}$, we have $\leb^{2(d-1)+5}(X') \gtrsim \leb^{2(d-1)+5}(X)$. Let $\tilde{\mu}^{2} = \lambda^{2-k_0\epsilon}$.

Since each $t_i \in S$, and each $|t_1 - t_2| \approxeq \tilde{\mu}^2$, we see that the $(t_0,t_\infty,t_1,t_{1'},t_2)$ reside in a set of measure $\lessapprox \leb^1(S)^4\tilde{\mu}^2$. Thus we may choose a $(t_0,t_\infty,t_1,t_{1'},t_2)$ so that, letting
\[
F'' = \{(\xi,x) : (\xi,x,t_0,t_\infty,t_1,t_{1'},t_2) \in X'\},
\]
we have 
\begin{equation} \label{F''estimate}
\leb^{2(d-1)}(F'') \gtrapprox (\leb^{1}(S))^{-4} \tilde{\mu}^{-2} \leb^{2(d-1)}(F) \lambda^{6}.
\end{equation}

We now define $G$ by the condition $G_X = F''$ and recall that
\[
|G| = \leb^{2(d-1)}(F''),\text{\ and\ } \ma(G) = \ma(F'').
\]
Let $t_{2'} = w_{t_0,t_\infty}(t_1,t_{1'},t_2)$. Since $s_{t_0,t_\infty}(t_1,t_{1'}) = s_{t_0,t_\infty}(t_2,t_{2'})$, we may apply Corollary \ref{sixsliceestimate} 
to obtain
\begin{align} \label{F''ub}
\leb^{2(d-1)}(F'')^{4} &\lesssim \left(D_{t_1,t_2} D_{t_0,t_\infty} D_{t_{1'},t_0}, D_{t_2',t_0} \right) \ma(F'')
\\ \notag &\ \ \ \ \ \ \ \cdot \sup_{t = t_{\infty}, t_0, t_1,t_{1'},t_2,t_{2'}} \leb^{d-1}(\pi_{t}(G))^7
\\ \notag &\lessapprox \tilde{\mu}^{-2(d-1)} \ma(F) \sup_{t = t_{\infty}, t_0, t_1,t_{1'},t_2,t_{2'}} \leb^{d-1}(\pi_{t}(G))^7.
\end{align}

By definition of $F''$, $\pi_{t_i}(G) \subset (E' \cap \gamma_{d}^{-1}(t_i))$
for $i \in \{0,\infty,1,1',2,2'\}$. Thus, we may combine (\ref{F''estimate}) and (\ref{F''ub}) to obtain
\begin{align}
\lefteqn{\leb^{2(d-1)}(F)^{4} \lambda^{24} \tilde{\mu}^{2(d-1)-8} \ma(F)^{-1}}\ \ \ \ \ \ \ \ \ \ \ \ \ \ \ \ \ \ \ \ \ \    
\\ \notag &\lessapprox \leb^{1}(S)^{16} \sup_{t = t_{\infty}, t_0, t_1,t_{1'},t_2,t_{2'}} \leb^{d-1}(E' \cap \gamma_{d}^{-1}(t))^7
\\ \notag &\lessapprox \leb^{d}(E)^{7} \leb^{2(d-1)}(F)^{-7\epsilon}.
\end{align}
Since $\tilde{\mu} \gtrsim \lambda$, we thus have
\[
\leb^{2(d-1)}(F)^{\frac{4}{7} + C\epsilon} \lambda^{\frac{14 + 2d}{7} + C\epsilon} \ma(F)^{-\frac{1}{7}} \lesssim \leb^{d}(E).
\]
Since $d \geq 6$, we have $4d+3 \geq 14 + 2d + C\epsilon$ and hence
\[
\left(\leb^{2(d-1)}(F)\right)^{\frac{4}{4d+3} + \epsilon} \lambda \ma(F)^{\frac{3}{4d+3} - \frac{4}{4d+3}} \lesssim \leb^d(E)^{\frac{7}{4d+3}}.
\]

\section{Further iteration} \label{pigeonholesection2}
By applying Proposition \ref{iterativemainestimate} once to 
(\ref{rpequal2}), we obtained Corollary \ref{sixsliceestimate}. 
One obtains the corollary below from $N-1$ iterative applications of Proposition \ref{iterativemainestimate}. This results in an improved value of $\frac{k}{\alpha - \beta}$, but also requires a larger collection of slices which satisfy a more complicated set of conditions.   
Recall the definition 
\[
u_{t_0,t_\infty}(t_i,\tilde{s}) = t_0 + \frac{(t_\infty - t_0)(t_i - t_0)}{\tilde{s}(t_i - t_\infty)},
\]
and note that $s_{t_0,t_\infty}(t_i,u_{t_0,t_\infty}(t_i,\tilde{s})) = \tilde{s},$  where $s_{t_0,t_\infty}$ is as defined in (\ref{iteratetildesdef}). 
\begin{corollary} \label{iterativecorollary}
Let $N \geq 3$ and $\sigma = (t_{0,1},t_{\infty,1},s_1,\ldots,t_{0,N-1},t_{\infty,N-1},s_{N-1},t_{0,N},t_{\infty,N}) \in \rea^{3N-1}$.
Let $\Gamma_{N+1}(\sigma) = \emptyset$, and for $1 \leq i \leq N$ let 
\[
\Gamma_{i}(\sigma) = \{t_{0,i},t_{\infty,i}\} \cup \Delta_i(\sigma) \cup \Gamma_{i+1}(\sigma)
\]
and
\[
\Delta_i(\sigma) = \{u_{t_{0,i},t_{\infty,i}}(t,s_{i}) : t \in \Gamma_{i+1}\}.
\] 
Suppose that 
\[
t_{0,i} \notin \Delta_i(\sigma) \cup \Gamma_{i+1}(\sigma) \cup \{t_{\infty,i}\}
\]
and 
\[
t_{\infty,i} \notin \Gamma_{i+1}(\sigma)
\]
for each $1 \leq i \leq N$.
Then for any set of lines $G$
\begin{equation} \label{iterativecorollaryconc}
|G|^{\alpha_N} \lesssim 
\left(\sup_{\tilde{t},\hat{t}} D_{\tilde{t},\hat{t}}^{\alpha_N}\right) 
\ma(G)^{\beta_N}
\left(\sup_{t} \leb^{d-1}(\pi_{t}(G))^{k_N}\right)
\end{equation}
where the right $\sup$ ranges over $t \in \Gamma_1(\sigma)$, where the left $\sup$ ranges over $\tilde{t},\hat{t}$ such that 
\[
\tilde{t} = t_{0,i},\ \hat{t} \in \{t_{\infty,i}\} \cup \Delta_i(\sigma),\ 1 \leq i \leq N,
\]
and where
\[
\alpha_{1} = 1,\ \ \beta_{1} = 0,\ \, k_{1} = 2, 
\]
and 
\[
\alpha_{i+1} = 2k_i,\ \ \beta_{i+1} = \beta_i + k_i - \alpha_i,\ \, \text{\ and\ } k_{i+1} = 4 k_{i} - \alpha_i
\]
for $1 \leq i < N.$
\end{corollary}

One may calculate the formulas  
\begin{gather*}
k_{i+1} = 4 k_i - 2 k_{i-1} \\
\beta_{i+1} = 2 k_{i-1}
\end{gather*}
which give 
\[
\frac{k_{i+1}}{\alpha_{i+1} - \beta_{i+1}} = 1 + \left(1 - \frac{k_{i-1}}{k_i}\right)^{-1}.
\]
Since $\frac{k_i}{k_{i-1}} = 4 - 2 \left(\frac{k_{i-1}}{k_{i-2}}\right)^{-1}$, the Banach contraction principle tells us that $\lim_{i \rightarrow \infty} \frac{k_{i}}{k_{i-1}} = 2 + \sqrt{2}.$
Thus
\begin{equation} \label{prfraclimit}
\lim_{i \rightarrow \infty} \frac{k_{i}}{\alpha_i - \beta_i} = 1 + \sqrt{2}. 
\end{equation}
Similarly
\begin{equation} \label{pqfraclimit}
\lim_{i \rightarrow \infty} \frac{k_{i}}{\alpha_i} = 1 + \frac{\sqrt{2}}{2}. 
\end{equation}

In the remainder of this section we use Corollary \ref{iterativecorollary} to show that the estimate (\ref{rwtestimate}) holds with $p_N,q_N,r_N$ satisfying
\[
\frac{r_N}{p_N} \geq \frac{k_N}{\alpha_N - \beta_N} - \epsilon,\ \text{\ and\ }
\frac{q_N}{p_N} \geq \frac{k_N}{\alpha_N} - \epsilon
\] 
where $\epsilon$ may be taken arbitrarily small. Thus, we obtain Theorem \ref{iteratedtheorem} from (\ref{prfraclimit}), (\ref{pqfraclimit}), and Claim \ref{rweclaim}. Due to the complicated nature of the set of slices $\Gamma_1(\sigma)$, we are not able to obtain an appropriately small value of $p_N$. Hence, for the sake of exposition, we will simplify the argument by passing up several opportunities to slightly improve $p_N$. For example, we will not employ the two-ends reduction. 

Let $E'$, $F'$, $\lambda$, $S$, and $S_{\xi,x}$ be as in Section \ref{pigeonholesection}. 
For $(\xi,x) \in F'$ let
\begin{align} \label{Xxixdef}
X_{\xi,x} &= \left\{\vphantom{\int}\right. \sigma \in \left(S_{\xi,x} \times S_{\xi,x} \times [-C \lambda^{-C_N},C \lambda^{-C_N}]\right)^{N-1} \times S_{\xi,x} \times S_{\xi,x}  
\\ \notag &\ \ \ \ \ \ \text{such\ that\ } \Gamma_{1}(\sigma) \subset S_{\xi,x},\\ \notag &\ \ \ \ \ \ \text{and\ such\ that\ }  |\tilde{t} - \hat{t}| \gtrsim \lambda^{C'_N} 
\\ \notag &\ \ \ \ \ \ \text{for all\ } \tilde{t} = t_{0,i},\ \hat{t} \in  \Delta_i(\sigma) \cup \Gamma_{i+1}(\sigma) \cup \{t_{\infty,i}\},
\\ \notag &\ \ \ \ \ \ \text{and for all\ } \tilde{t} = t_{\infty,i}, \hat{t} \in \Gamma_{i+1}(\sigma), 
\\ \notag &\ \ \ \ \ \ \text{where\ }1 \leq i \leq N \left. \vphantom{\int}\right\},
\end{align}
where $C_N$ and $C'_N$ will be determined below.
One should think of $X_{\xi,x}$ as the set of candidates for $\sigma$
 in Corollary \ref{iterativecorollary}.
Our aim is to find a lower bound for $\leb^{3N-1}(X_{\xi,x})$. This will be accomplished by providing the lower bound for a subset $Y_N$ of $X_{\xi,x}$ which is appropriately compatible with the following estimate. 

\begin{claim} \label{unifclaim}
Let $I \subset [0,1]$ with $\leb^1(I) < \infty$.
For $(t_0,t_\infty,s) \in I \times I \times \rea$
let 
\begin{multline} \label{unifclaimeq1}
\widetilde{I}_{t_0,t_\infty,s} = \{t \in I : u_{t_0,t_\infty}(t,s) \in I, \text{\ and\ } |\tilde{t} - \hat{t}| \gtrsim \leb^1(I) 
\\ \text{\ for\ } \tilde{t} \in \{t, u_{t_0,t_\infty}(t,s) \}, \hat{t} \in \{t_0,t_\infty\} \}.
\end{multline}
Then letting
\begin{multline} \label{unifclaimeq2}
\mathcal{P}(I) = \{(t_{0},t_\infty,s) \in I \times I \times [-C \leb^1(I)^{-2}, C \leb^1(I)^{-2}] : 
\\ |t_0 - t_\infty| \gtrsim \leb^1(I), \text{\ and\ } \leb^{1}(\tilde{I}_{t_0,t_\infty,s}) \gtrsim \leb^1(I)^{6}
  \}, 
\end{multline}
we have
\begin{equation} \label{unifclaimeq3}
\int_{\mathcal{P}(I)} \leb^{1}(\widetilde{I}_y)\ dy \gtrsim \leb^1(I)^6.
\end{equation}
\end{claim}

\begin{proof}
We argue as in Section \ref{pigeonholesection}.
Let 
\[
P = \{(t_0,t_\infty) \in I \times I : |t_0 - t_\infty| \gtrsim \leb^1(I) \}
\]
and note that
\begin{equation} \label{measP}
\leb^2(P) \gtrsim \leb^{1}(I)^2.
\end{equation} 
For each $(t_0,t_\infty) \in P$, we let 
\[
Q_{t_0,t_\infty} = \{(t_1,t_{1'}) \in I \times I
: |t_i - t_j| \gtrsim \leb^1(I) \ \text{for\ all\ } i \neq j \in \{0,\infty,1,1'\}\}.
\] 
and note that 
\begin{equation} \label{measQ}
\leb^{2}(Q_{t_0,t_\infty}) \gtrsim \leb^1(I)^2.
\end{equation}
Consider any fixed $t_0,t_\infty \in P$, and let $Q = Q_{t_0,t_\infty}$. Then, as in (\ref{changevarQequation}), we have
\begin{align} \label{changevarQ2} 
\leb^{2}(Q) &= \int_{\rea} \int_{\rea} \chi_{Q}(t_{1},t_{1'})\ dt_1\ dt_{1'} 
\\ \notag &= \int_{\rea} \int_{\rea} \chi_{Q}\left(t_1,u_{t_0,t_\infty}(t_1,s')\right)\left|\frac{(u_{t_0,t_\infty}(t_1,s') - t_0)^2(t_1 - t_\infty)}{(t_1 - t_0)(t_\infty - t_0)} \right|\ dt_1\ ds'
\\ \notag &\lesssim \left(\leb^{1}(I)\right)^{-2}
 \int_{\rea} \int_{\rea} \chi_{Q}\left(t_1,u_{t_0,t_\infty}(t_1,s')\right)\ dt_1\ ds'.
\end{align}
Note that in $Q$ each $|s'| \lesssim \leb^1(I)^{-2}$.
Thus, letting 
\[
\overline{S}_{t_0,t_\infty} = \left\{s' : \int_{\rea} \chi_{Q}\left(t_1,u_{t_0,t_\infty}(t_1,s')\right)\ dt_1 \ll \leb^1(I)^6 \right\}
\]
we have  
\begin{equation} \label{negligibleS}
\int_{\overline{S}_{t_0,t_\infty}} \int_{\rea} \chi_{Q}\left(t_1,u_{t_0,t_\infty}(t_1,s')\right)\ dt_1\ ds' \ll  
\leb^{1}(I)^4 \lesssim \leb^{2}(Q)\leb^1(I)^2.
\end{equation}
Next, we note that if $s \notin \overline{S}_{t_0,t_\infty}$ then 
$(t_0,t_\infty,s) \in \mathcal{P}(I)$.
Thus, 
\begin{align*}
\int_{\mathcal{P}(I)} \leb^{1}(\tilde{I}_{x}) dx &\geq \int_{P} \int_{\rea \setminus \overline{S}_{t_0,t_\infty}} \left(\int_{\rea} 
\chi_{Q_{t_0,t_\infty}}\left(t_1,u_{t_0,t_\infty}(t_1,s')\right)\ dt_1 \right)\ ds'\ dt_{0}\ dt_{\infty}
\\ &\gtrsim \leb^{1}(I)^6,
\end{align*}
where the second inequality follows from (\ref{measP}), (\ref{measQ}), (\ref{changevarQ2}), and (\ref{negligibleS}).
\end{proof}

Taking $I = S_{\xi,x}$ in Claim \ref{unifclaim}, we let $Y_1 = \mathcal{P}(S_{\xi,x})$ and for each $(t_0,t_\infty,s) = y \in Y_1$ let $I_y = \widetilde{(S_{\xi,x})}_y$. Recalling that $\leb^1(S_{\xi,x}) \gtrsim \lambda$, we 
have by definition of $\mathcal{P}(S_{\xi,x})$
\begin{equation} \label{yinductionstart2}
\leb^1(I_y) \gtrsim \lambda^6 \text{\ for\ } y \in Y_1
\end{equation}
and we see from (\ref{unifclaimeq3}) that
\begin{equation} \label{yinductionstart1}
\int_{Y_1} \leb^1(I_y)\ dy \gtrsim \lambda^6.
\end{equation}
For $j = 2,\ldots,N-1$, we define $Y_j$ and $I_y$ recursively, letting
\[
Y_{j} = \{(y',y'') : y' \in Y_{j-1} \text{\ and\ } y'' \in \mathcal{P}(I_{y'})\} \subset \rea^{3j}
\]
and 
\[
I_{(y',y'')} = \widetilde{(I_{y'})}_{y''} \text{\ for\ } (y',y'') \in Y_{j}. 
\]
From (\ref{yinductionstart2}), the definition of $\mathcal{P}(I),$ and induction, we see that 
\begin{equation} \label{measIyeq}
\leb^{1}(I_y) \gtrsim \lambda^{6^j} \text{\ for \ } y \in Y_{j}.
\end{equation}
From (\ref{unifclaimeq3}), (\ref{yinductionstart1}), (\ref{measIyeq}), and induction, we have
\begin{align} \label{measYjeq}
\int_{Y_j} \leb^1(I_y)\ dy &= \int_{Y_{j-1}} \int_{\mathcal{P}(I_{y'})} \leb^{1}(\widetilde{(I_{y'})}_{y''})\ dy''\ dy'
\\ \notag &\gtrsim \int_{Y_{j-1}} (\leb^{1}(I_{y'}))^{6}\ dy'
\\ \notag &\gtrsim \lambda^{5 \cdot 6^{j-1}} \int_{Y_{j-1}} \leb^{1}(I_{y'}) \ dy'
\\ \notag &\gtrsim \lambda^{5 \cdot 6^{j-1}} \lambda^{6^{j-1}} = \lambda^{6^j}.
\end{align}
Finally, we let
\[
Y_N = \{(y,t,t') : y \in Y_{N-1}; t,t' \in I_y; \text{\ and\ } |t - t'| \gtrsim \leb^1(I_y)\}.
\]
From (\ref{measIyeq}) and (\ref{measYjeq}), we have
\[
\leb^{3N - 1}(Y_N) \gtrsim \lambda^{2 \cdot 6^{N-1}}.
\]

We will now verify that $Y_N \subset X_{\xi,x}$, where in the definition (\ref{Xxixdef}) of 
$X_{\xi,x}$ we have $C_N = 2 \cdot 6^{N-2}$ and $C'_N = 6^{N-1}$. 
Let
\[
\sigma_N =
(t_{0,1},t_{\infty,1},s_1,\ldots,t_{0,N-1},t_{\infty,N-1},s_{N-1},t_{0,N},t_{\infty,N}) \in Y_N
\]
and for $j = 1,\ldots,N-1$ let
\[
\sigma_j = (t_{0,1},t_{\infty,1},s_1,\ldots,t_{0,j},t_{\infty,j},s_{j}),
\]
where we note that $\sigma_j \in Y_j$.
Additionally, define $I_{\sigma_0} = S_{\xi,x}.$
By the definition of $Y_i$, we have each  
\begin{equation} \label{tiinpi}
(t_{0,i},t_{\infty,i},s_i) \in \mathcal{P}(I_{\sigma_{i-1}}).
\end{equation}
Since, by (\ref{measIyeq}), each $\leb^1(I_{\sigma_{i-1}})^{-2} \lesssim \lambda^{-C_{N}}$, we thus have, by (\ref{unifclaimeq2}),
\[
\sigma_N \in (S_{\xi,x} \times S_{\xi,x} \times [-C \lambda^{-C_N},C \lambda^{-C_N}])^{N-1} \times S_{\xi,x} \times S_{\xi,x}.
\]
Next, we note that $\Gamma_N(\sigma_N) = \{t_{0,N},t_{\infty,N}\} \subset I_{\sigma_{N-1}}$ and that, since $I_{\sigma_{i-1}} = \widetilde{(I_{\sigma_{i-2}})}_{t_{0,i-1},t_{\infty,i-1},s_{i-1}}$,
\[
\Gamma_{i}(\sigma_N) \subset I_{\sigma_{i-1}} \Rightarrow \Delta_{i-1}(\sigma_N) \subset I_{\sigma_{i-2}} \text{\ and\ thus\ } \Gamma_{i-1}(\sigma_N) \subset I_{\sigma_{i-2}}.
\]
So by induction $\Gamma_{i}(\sigma_N) \subset I_{\sigma_{i-1}}$ for $1 \leq i \leq N$, and in particular $\Gamma_1(\sigma_N) \subset I_{\sigma_0} = S_{\xi,x}$.
Again using (\ref{tiinpi}) and (\ref{unifclaimeq2}), we obtain
\[
|t_{0,i} - t_{\infty,i}| \gtrsim \leb^1(I_{\sigma_{i-1}})
\gtrsim \lambda^{C'_N}.
\]
Also, since $\Gamma_{i+1}(\sigma_N) \subset I_{\sigma_{i}} = \widetilde{(I_{\sigma_{i-1}})}_{t_{0,i},t_{\infty,i},s_i}$, we see that
\[
|\tilde{t} - \hat{t}| \gtrsim \leb^1(I_{\sigma_{i-1}})
\gtrsim \lambda^{C'_N} \text{\ for\ } \tilde{t} \in \{t_{0,i},t_{\infty,i}\},\ \hat{t} \in \Delta_{i}(\sigma_N) \cup \Gamma_{i+1}(\sigma_N).
\]
Thus $\sigma_N \in X_{\xi,x}$ and $Y_N \subset X_{\xi,x}$. In particular 
\[
\leb^{3N-1}(X_{\xi,x}) \gtrsim \lambda^{2\cdot6^{N-1}}.
\]
 
Let
\[
X = \{(\xi,x,\sigma) : (\xi,x) \in F' \text{\ and\ } \sigma \in X_{\xi,x} \}
\]
and note that 
\[
\leb^{2(d-1) + 3N - 1}(X) \gtrsim \lambda^{2\cdot6^{N-1}} \leb^{2(d-1)}(F').
\]
Since the $\sigma$'s reside in a set of measure $\lesssim \leb^1(S)^{2N} \lambda^{-(N-1)C_N}$, 
we see that, letting $C''_N = 2\cdot6^{N-1} + (N-1)C_N$, we may find 
a fixed $\sigma$ so that 
\begin{equation} \label{iterativeFeq}
\leb^{2(d-1)}(F'') \gtrsim \lambda^{C''_N} \leb^1(S)^{-2N} \leb^{2(d-1)}(F')
\end{equation}
where
\[
F'' = \{(\xi,x) : (\xi,x,\sigma) \in X \}.
\]
Since $\mass(E,F') \gtrsim \lambda \leb^{2(d-1)}(F)$ and $T[\chi_E] \leq 1$, 
we have 
\begin{equation} \label{iterativeFeq2}
\leb^{2(d-1)}(F') \gtrsim \lambda {\leb^{2(d-1)}}(F).
\end{equation}

We now apply Corollary \ref{iterativecorollary} with the set of lines 
$G$ defined by $G_X = F''$. We then have
\begin{gather*}
|G| = \leb^{2(d-1)}(F'')
\intertext{and}
\ma(G) = \ma(F'') \leq \ma(F).
\end{gather*} 
Also, since $\Gamma_{1}(\sigma) \in S_{\xi,x}$ for 
every $(\xi,x) \in F''$, we have $\pi_t(G) \subset E \cap \gamma_{d}^{-1}(S)$ for every $t \in \Gamma_1(\sigma)$. Furthermore, by definition of $X_{\xi,x}$, we have 
\[
\sup_{\tilde{t},\hat{t}} D_{\tilde{t},\hat{t}} \lesssim \lambda^{-(d-1)C'_N}.
\]
where the $\sup$ ranges over
\[
\tilde{t} = t_{0,i},\ \hat{t} \in \{t_{\infty,i}\} \cup \Delta_i(\sigma),\ 1 \leq i \leq N.
\]
Thus, from (\ref{iterativecorollaryconc}) we obtain
\begin{multline} \label{iterativecorollaryagain}
\leb^{2(d-1)}(F'')^{\alpha_N} \lesssim 
\\ \left( \lambda^{-(d-1)C'_N}\right)^{\alpha_N} 
\ma(F)^{\beta_N}
\left(\sup_{t \in S} \leb^{d-1}(E \cap \gamma_d^{-1}(t))^{k_N}\right).
\end{multline}
Noting that $\frac{k_N}{\alpha_N} < 2$ (and certainly $< 2N$)
and
\[
\leb^1(S) \left(\sup_{t \in S} \leb^{d-1}(E) \cap \gamma_d^{-1}(t)\right) \lesssim \lambda^{-\epsilon} \leb^{2(d-1)}(F)^{-\epsilon} \leb^{d}(E),
\]
we obtain from 
(\ref{iterativeFeq}), (\ref{iterativeFeq2}), and (\ref{iterativecorollaryagain})
\[ 
\leb^{2(d-1)}(F)^{\alpha_N + k_N \epsilon} \lambda^{C'''_N \alpha_N + k_N \epsilon} 
\ma(F)^{- \beta_N}
\lesssim
\leb^{d}(E)^{k_N}.
\]
where $C'''_N =1+(d-1)C'_N + C''_N$.

Finally, this gives 
\[
\lambda \leb^{2(d-1)}(F)^{\frac{1}{q_N}} \ma(F)^{\left(\frac{1}{q_N} - \frac{1}{r_N}\right)}
\lesssim \leb^d(E)^{\frac{1}{p_N}}
\]
where
\begin{gather}
p_N = \frac{C'''_N\alpha_N + k_N\epsilon}{k_N} \geq 6^{N-2}(6(d+1) + 2(N-1)) \frac{\alpha_N}{k_N} + \epsilon
\\ \notag q_N = \frac{C'''_N \alpha_N + k_N\epsilon}{\alpha_N + k_N \epsilon} \geq \left(\frac{k_N}{\alpha_N} - C \epsilon\right)p_N  
\\ \notag r_N = \frac{C'''_N \alpha_N + k_N \epsilon}{\alpha_N - \beta_N + k_N \epsilon} \geq \left(\frac{k_N}{\alpha_N - \beta_N} - C\epsilon\right)p_N.  
\end{gather}
Thus, from (\ref{prfraclimit}), we see that by taking $N$ large and $\epsilon$ small, we have $\frac{r}{p}$ arbitrarily close to $1+\sqrt{2}$ (although this comes with the price of a very large $p$).

\section{Nikodym mixed-norms} \label{projectivesection}
We now follow a method of Tao from \cite{tao} to show that Corollary \ref{nikodymcorollary} follows from Theorem \ref{maintheorem}.
For $z \in \rea^d$ write $x_z = \proj_H(z)$ and $t_z e_d = \proj_{e_d}(z)$
where $\proj\ $ denotes orthogonal projection.
For nonnegative integers $j$ let 
\[
S_j = \{z \in \rea^d : 2^{-(j+1)} < t_z \leq 2^{-j}\}.
\] 
We first prove 
Corollary \ref{nikodymcorollary} in the special case when $f$ is supported 
on $S_0$. In fact, to prove this case it suffices, since $T$ is local and 
$p \leq q \leq r$, to consider the case when $f$ is supported on $S_0 \cap 
Q$ where $Q$ is the cube centered at $\frac{1}{2}e_d$ with side length $1$.
Furthermore, assume that $f$ is positive.

We then consider the projective
transformation
\[
\phi(z) = \frac{x_z + e_d}{t_z}.
\]
The idea is that we have $T[f](\xi,x) \approx T[f \composed \phi](x,\xi)$.
To line everything up, we will in fact estimate 
$T[f \composed \phi \composed d_2]$ where $d_2(y) = 2y$. Then
\begin{align*}
T[f \composed \phi \composed d_2](\xi,x) &= \int_{\frac{1}{2}}^1 f \composed \phi \composed d_2 (x + t(\xi + e_d))\ dt
\\ &= \int_{\frac{1}{2}}^1 f\left(\xi + \frac{1}{2t}\left(2x + e_d\right) \right)\ dt
\\ &\approx \int_{\frac{1}{2}}^1 f(\xi + \tilde{t}(2x + e_d))\ d\tilde{t}
\\ &= T[f](2x, \xi),
\end{align*}
where the first and last equations follow from the fact that $f$ is supported on $S_0$.
Thus by Theorem \ref{maintheorem}
\begin{align*}
\lefteqn{\left(\int_{\rea^{d-1}} \left( \int_{B(0,C)} |T[f](\xi,x)|^r\ d\xi \right)^{\frac{q}{r}}\ dx \right)^{\frac{1}{q}}}
\\ &=\left(\int_{B(0,C')} \left( \int_{B(0,C)} |T[f](\xi,x)|^r\ d\xi \right)^{\frac{q}{r}}\ dx \right)^{\frac{1}{q}}
\\ &\lesssim \left(\int_{B(0,C')} \left( \int_{B(0,C)} \left|T[f \composed \phi \composed d_2]\left(x,\frac{1}{2}\xi\right)\right|^r\ d\xi \right)^{\frac{q}{r}}\ dx \right)^{\frac{1}{q}}
\\ &\lesssim \left(\int_{B(0,C')} \left( \int_{\rea^{d-1}} |T[f \composed \phi \composed d_2](x,\xi)|^r\ d\xi \right)^{\frac{q}{r}}\ dx \right)^{\frac{1}{q}}
\\ &\lesssim \|f \composed \phi \composed d_2\|_{L^p(\rea^d)}.
\end{align*}
where we use the fact that $f$ is supported on $Q$ for the first equation.
Since $f$ is supported on $S_0$, we have
\[
\|f \composed \phi \composed d_2\|_{L^p(\rea^d)} \lesssim \|f\|_{L^p(\rea^d)}.
\]

Lifting our support assumptions on $f$, we note that
\[
\|T[f]\|_{L^q(L^r),N} \leq \sum_{j=0}^\infty \|T[\chi_{S_j} f]\|_{L^q(L^r),N}.
\]
However each $\chi_{S_j} f \composed d_{2^{-j}}$ is supported on $S_0$,
and  
\[
T[\chi_{S_j}f](\xi,x) = 2^{-j} T[\chi_{S_j}f\composed d_{2^-j}](\xi,2^j x).
\]
Thus, for each $j$
\begin{align*}
\|T[\chi_{S_j} f]\|_{L^q(L^r),N} &= 2^{-j(1 + \frac{d-1}{q})}\|T[\chi_{S_j} f \composed d_{2^{-j}}]\|_{L^q(L^r),N}
\\ &\lesssim  2^{-j(1 + \frac{d-1}{q})} \|\chi_{S_j} f \composed d_{2^{-j}}\|_{L^{p}(\rea^d)}
\\ &\leq 2^{-j(1 + \frac{d-1}{q} - \frac{d}{p})} \|f\|_{L^{p}(\rea^d)}.
\end{align*}
Hence, provided (\ref{qpcondition}) holds with strict inequality, which is 
indeed the case here, we have
\[
\sum_{j=0}^\infty \|T[\chi_{S_j} f]\|_{L^q(L^r),N} \lesssim \|f\|_{L^{p}(\rea^d)}.
\]

\end{document}